\documentclass[11pt]{amsart}

\usepackage{color}
\usepackage{amsmath}
\usepackage{amsxtra}
\usepackage{amscd}
\usepackage{amsthm}
\usepackage{amsfonts}
\usepackage{amssymb}
\usepackage{eucal}
\usepackage{epsfig}
\usepackage{graphics}
\usepackage[matrix,arrow]{xy}

\usepackage{hyperref}

\theoremstyle{plain}
\newtheorem{theorem}{Theorem}[section]

\newtheorem{prop}[theorem]{Proposition}

\newtheorem{conj}[theorem]{Conjecture}
\newtheorem{defi}[theorem]{Definition}
\newtheorem{example}[theorem]{Example}

\newcommand{\g}{\mathfrak{g}}

\newcommand{\Psib}{\mbox{\boldmath$\Psi$}}

\newcommand{\nc}{\newcommand}
\nc{\on}{\operatorname}
\nc{\la}{\lambda}
\nc{\wh}{\widehat}
\nc{\wt}{\widetilde}
\nc{\sw}{{\mathfrak s}{\mathfrak l}}
\nc{\ghat}{\wh{\g}}
\nc{\hhat}{\wh{\h}}
\nc{\mc}{\mathcal}
\nc{\bi}{\bibitem}
\nc{\pa}{\partial}
\nc{\ppart}{(\!(t)\!)}
\nc{\pparl}{(\!(\la)\!)}
\nc{\zpart}{(\!(z^{-1})\!)}
\nc{\n}{{\mathfrak n}}
\nc{\ol}{\overline}
\nc{\mb}{\mathbf}
\nc{\bb}{{\mathfrak b}}
\nc{\su}{\wh\sw_2}
\nc{\h}{{\mathfrak h}}
\nc{\can}{\on{can}}
\nc{\ntil}{\wt{\n}}
\nc{\pone}{{\mathbb P}^1}
\nc{\bs}{\backslash}
\nc{\al}{\alpha}
\nc{\gt}{{\mathfrak g}'}
\nc{\ds}{\displaystyle}

\nc{\Sig}{{\mb \Sigma}}

\theoremstyle{definition}
\newtheorem{rem}[theorem]{Remark}

\begin{document}

\title{Symmetries of Grothendieck rings in representation theory}

\author[David Hernandez]{David Hernandez}

\address{IMJ-PRG, Universit\'e Paris Cit\'e and Sorbonne Universit\'e, CNRS, F-75006, Paris, France}

\begin{abstract}  This is a written version of the invited lecture at the 9th European Congress of Mathematics in July 2024 in Sevilla.
We review certain new symmetries of Grothendieck rings that have emerged in representation theory.
\end{abstract}

\maketitle

\tableofcontents

Categories of representations of groups, or more generally, bialgebras and quantum groups, are endowed with monoidal structures. The Grothendieck ring of such a category has a basis consisting of classes of simple objects, with structure constants defined as multiplicities of simple representations in the decomposition of tensor products. The resulting rings together with their bases can have a very intricate structure, and many of them are fundamental with important applications 
(for example, leading to algorithms to compute the dimension of simple representations).

In recent years, new symmetries of Grothendieck rings have emerged, which can be used to better understand these rings. We will discuss these symmetries and illustrate some of their applications. 

First, cluster symmetries give powerful tools to investigate the structure of the Grothendieck rings \cite{HL1}. As an application, we have obtained the proof of long-standing conjectures on character formulas for simple finite-dimensional representations of quantum affine algebras (which are infinite-dimensional quantum groups) \cite{FHOO}. 

Secondly, motivated by symmetries of quantum integrable models arising in physics, we have established remarkable relations in associated Grothendieck rings \cite{FH}, generalizing the celebrated Baxter relations; these relations led us to the discovery of a surprising new action of the Weyl group on these Grothendieck rings \cite{FH2}. These Weyl group symmetries are also closely related to the cluster symmetries \cite{GHL}.

\newpage

\section{Grothendieck group}

To start with, let us remind the general construction Grothendieck groups.

\subsection{Grothendieck group of an additive category}

Let $\mathcal{C}$ be an additive category with distinguished short exact sequences
\begin{equation}\label{ese}A\rightarrow B \rightarrow C.\end{equation}
For example, one may consider the case of an abelian category with all ordinary exact sequences.

\begin{defi} The Grothendieck group $K(\mathcal{C})$ of $\mathcal{C}$ is the group generated by the elements $[M]$ for each isomorphism class of objects $M$ in $\mathcal{C}$, with the relations 
$$[A] + [C] = [B]$$
for each exact sequence (\ref{ese}).
\end{defi}

We obtain the important algebraic structure $K(\mathcal{C})$ attached to the category $\mathcal{C}$.
We will discuss many examples in this lecture.

The main original motivation appeared in algebraic geometry. Consider the $\mathcal{C}$ category of coherent sheaves 
on an algebraic variety $X$. The Grothendieck group 
$$K(X) = K(\mathcal{C})$$ 
is one of the starting points of $K$-theory in algebraic geometry.
Historically, such constructions were first studied by Grothendieck for categories of coherent and locally free sheaves on schemes, in proving the Riemann–Roch theorem.

Another basic but fundamental example is obtained for a field $k$ : the category $\mathcal{C} = \mathcal{V}$ of finite-dimensional $k$-vector spaces. The isomorphism classes are parametrized by the dimension. Then $\mathcal{V}$ is semi-simple and the direct sum descends to the sum in the Grothendieck group, that is for $V$ and $W$ in $\mathcal{V}$ we have
$$[V\oplus W] = [V] + [W].$$
This implies that $K(\mathcal{V})$ is freely generated by the class $[k]$ of the one-dimensional vector space 
$$K(\mathcal{V}) = \mathbb{Z} [k] \simeq \mathbb{Z} $$ 
as we have
$$[V] = \text{dim}(V) . [k]$$
for $V$ in $\mathcal{V}$.

\subsection{Grothendieck groups in representation theory}

We now discuss the Grothendieck group of finite-dimensional representations of an algebra.

Let $\mathcal{A}$ be an algebra over a field $k$ and $\mathcal{F}$ be the category of finite-dimensional representations of $\mathcal{A}$. 

We have the following classical and fundamental result.

\begin{theorem}\label{jhs} [Jordan-H\"older series]
Each object $M$ in $\mathcal{F}$ admits a series of subobjects
$$M = M_0\supset M_1\supset\cdots \supset M_N = \{0\}$$
with $M_i/M_{i+1}$ simple object. 

The number $n_V(M)$ of occurrence of a simple object $V$ (as a quotient $M_i/M_{i+1}$) depends only on $M$, and not on the series of 
subobjects. 
\end{theorem}

Such a series is called a composition series or a Jordan-H\"older series of the representation. A simple quotient $M_i/M_{i+1}$ is called a composition factor of $M$, and its number of occurrence is called its multiplicity. 

For completeness, let us briefly remind the arguments of the proof. 

\begin{proof} The existence of the series is obtained by induction on the dimension of the representation. 
The uniqueness of $n_V(M)$ is also proved by induction on the dimension : it is obvious for a simple representation. 
In general, assume that we have two series $M_i$ and $M_i'$ as above. If $M_1 = M_1'$, we can conclude using the induction hypothesis 
for this representation. Otherwise 
$$M\supset M_1 \supset M_1\cap M_1'$$ 
and 
$$M\supset M_1' \supset M_1\cap M_1'$$ 
are the first few terms of composition series with the same composition factors. Using the induction hypothesis, they have the same composition factors.
\end{proof}

One can then define the length of a module $M$ as the sum 
$$l(M) = \sum_{V \text{ simple class in }\mathcal{F}} n_V(M).$$

As a consequence of Theorem \ref{jhs}, we obtain that 
$$K(\mathcal{F}) = \bigoplus_{V \text{ simple class}} \mathbb{Z}[V]$$
has a canonical basis parametrized by simple classes. The Grothendieck group is a group with a distinguished basis. 

The coefficient on the basis are the multiplies of simple representations, that is for $M$ object in $\mathcal{F}$, one has
$$[M] = \sum_{V \text{ simple class}}n_V(M)[V].$$
In general, an element 
$$\sum_{V \text{ simple class}}n_V [V]$$
with all $n_V\geq 0$ is called the class of an {\it actual} representation.

If one of the $n_V < 0$, this is the class of a {\it virtual} representation. The Grothendieck group allows to 
treat in a uniform way actual and virtual representations.

\section{Monoidal categories}

We now assume that an abelian category $\mathcal{C}$ has a monoidal structure, that is we have 
tensor products $M\otimes N$ for $M,N$ objects in $\mathcal{C}$ with isomorphisms
$$(M\otimes N)\otimes P \simeq M\otimes (N\otimes P)$$
as well as a unit object ${\bf 1}$ in $\mathcal{C}$ so that for $M$ object in $\mathcal{C}$ 
$$M\otimes {\bf 1}\simeq \mathbb{\bf 1}\otimes M\simeq M.$$
We have additional compatibility properties with the abelian structure (see \cite{ML}).

\begin{rem} We do not assume that the category is necessarily a tensor category, that is we do not assume the category admits dualities compatible with the monoidal structure.\end{rem}

Let us discuss important examples of monoidal categories.

\medskip

1) The category $\mathcal{V}$ of finite-dimensional $k$-vector spaces has a natural monoidal structure.

\medskip

2) For a group $G$, the category of finite-dimensional representations of $G$ (that is of its group-algebra $k[G]$) has 
a monoidal structure. The action of $g\in G$ on a pure tensor $v\otimes w\in V\otimes W$ in a tensor product of representations 
is defined by 
$$g.(v\otimes w) = (g.v)\otimes (g.w).$$

\medskip

3) For a Lie algebra $\mathfrak{g}$, the category of finite-dimensional representations of $\mathfrak{g}$ (that is of 
its universal enveloping algebra $\mathcal{U}(\mathfrak{g})$) has a monoidal structure. 
The action of $g\in \mathfrak{g}$ on a pure tensor $v\otimes w\in V\otimes W$ in a tensor product of representations 
is defined by 
$$g.(v\otimes w) = (g.v)\otimes w + v\otimes (g.w).$$
This is an infinitesimal version of the group action above.

\medskip

4) More generally : for a Hopf algebra $H$, that is an algebra with a coproduct and a counit
$$\Delta : H\rightarrow H\otimes H\text{ and }\epsilon : H\rightarrow k.$$
The category of finite-dimensional representations of $H$ has a monoidal structure. This includes all examples above, 
and new examples, in particular the Drinfeld-Jimbo {\bf quantum groups} $$\mathcal{U}_q(\mathfrak{g})$$
associated to $q\in \mathbb{C}^*\setminus\{-1, 1\}$ quantum parameter and 
$\mathfrak{g}$ simple Lie algebra. 
It is a $q$-deformation of $\mathcal{U}(\mathfrak{g})$, originated from mathematical physics. For example, $\mathcal{U}_q(\mathfrak{sl}_2)$ 
is the algebra by generators $E$, $F$, $K^{\pm 1}$ and relations 
$$KE = q^2EK\text{ , }KF = q^{-2}FK\text{ and }[E,F] = \frac{K - K^{-1}}{q - q^{-1}}.$$

\medskip

5) Hecke algebras, of their variations such as quiver Hecke algebras (as defined by Khovanov-Lauda and Rouquier) 
quiver rise categories of representations which have a monoidal structures via convolutions.

\medskip

6) In various geometric contexts, convolution diagrams 
$$X_1\overset{p}{\leftarrow} X_2\overset{q}{\rightarrow} X_3$$
of regular morphisms $p$, $q$ between remarkable varieties allow to consider convolution 
products and then monoidal structures on certain categories of sheaves (see for instance \cite{CG, BR} and 
references therein).

\section{Grothendieck ring}

\subsection{Grothendieck ring of a monoidal category}

Let us define the Grothendieck ring of $\mathcal{C}$ monoidal category. 

\begin{prop} The Grothendieck group $K(\mathcal{C})$ of $\mathcal{C}$ inherits a ring structure such that 
$$[M].[N] = [M\otimes N]$$
for $M$, $N$ objects in $\mathcal{C}$ and the class $[{\bf 1}]\in K(\mathcal{C})$ is the neutral element.
\end{prop}

\begin{example} $\mathcal{C} = \mathcal{V}$ the category of finite-dimensional vector spaces on a field $k$, then the group isomorphism
$$K(\mathcal{C})\simeq \mathbb{Z}$$ 
is a ring isomorphism. This follows for example from the formula $\text{dim}(V\otimes W) = \text{dim}(V)\text{dim}(W)$.
\end{example}

\begin{rem} If the category is braided, that is 
$$M\otimes N\simeq N\otimes M$$ 
for all objects $M,N$ in $\mathcal{C}$, then the Grothendieck ring is commutative. The converse statement is false : there 
are non braided monoidal categories whose Grothendieck ring is commutative. We will see examples below. 
\end{rem}

\subsection{Grothendieck rings in representation theory}

Let $\mathcal{C}$ be the category of finite-dimensional representations of a bialgebra algebra $H$. The Grothendieck ring 
$$K(\mathcal{C}) = \bigoplus_{V\text{ simple}}\mathbb{Z}[V]$$
has a canonical basis $([V])_{V\text{ simple}}$, and the constant structures on this basis are positive constant structures  : 
they are given by the multiplicities of simple modules in tensor products (sometimes called Clebsch-Gordan coefficients 
or fusion coefficients).

\begin{rem} If $H$ is a Hopf algebra, $\mathcal{C}$ is a tensor category, and dualities induce ring automorphisms of 
$K(\mathcal{C})$, which preserve the canonical basis.\end{rem}

\begin{example}\label{bex} let $\mathcal{C}$ be the category of finite-dimensional representations of the Lie algebra $\mathfrak{sl}_2(\mathbb{C})$. 
Then we have a ring isomorphism
$$K(\mathcal{C})\simeq \mathbb{Z}[X]$$
with the polynomial ring in one variable $X$. Here
$$X = [V]\text{ where $V = \mathbb{C}^2$ is the natural representation.}$$
It is well-known that $\mathfrak{sl}_2(\mathbb{C})$ has a simple complex representation of dimension $n + 1$ for any $n\geq 0$, and that 
this representation is unique up to isomorphism. Let us denote by $Q_n$ the isomorphism class of this representation in $K(\mathcal{C})$. Then 
the basis $(Q_n)_{n\geq 0}$ of simple classes satisfies the induction relation
\begin{equation}\label{qq}Q_n^2 = Q_{n+1}Q_{n-1} + 1.\end{equation} 
One can compute the canonical basis of simple classes 
$$(Q_0,Q_1,Q_2,\cdots) = (1, X, X^2 - 1, X^3 - 2X, \cdots),$$
which is different that the standard monomial basis of $\mathbb{Z}[X]$.
\end{example}

\begin{rem} Note that, independently on Grothendieck rings, it is not a priori completely that the 
inductive sequence (\ref{qq}) with initial conditions $Q_0 = 1$ and $Q_1 = X$, a priori in $\mathbb{Q}(X)$, defines a sequence 
in the ring $\mathbb{Z}[X]$.
\end{rem}

\subsection{Applications of Grothendieck rings in representation theory}

The study of study of Grothendieck rings of categories of representations lead to interesting developments in various directions. 
Let us list some of them.

\medskip

 {\it Representation Theory} : there is close relation between the understanding of the global structure of 
Grothendieck rings and the description of the structure of simple representations. 
For example, if one knows the structure of certain simple representations, one may extract informations on 
the simple constituents of their tensor products. 

\begin{example} This approach will be discussed below to compute the characters of simple representations of quantum affine algebras 
from the characters of fundamental representations, following \cite{Nak:quiver} and then \cite{H1, FHOO}.
\end{example}

\medskip

 {\it Categorification} : one can realize sometimes an algebra $\mathcal{A}$ as a Grothendieck ring:
$$\mathcal{A}\simeq K(\mathcal{C}).$$
The categorified structure $\mathcal{A}$ inherits a canonical basis from the canonical basis of the Grothendieck ring.

\begin{example} The coordinate ring $\mathbb{C}[N]$, with $N$ unipotent subgroup of a Lie group $G$, can be realized from representations of quiver Hecke algebras  \cite{KL, R}
or of quantum affine algebras \cite{HLCr}, leading to a categorical interpretation of its dual canonical basis. A basic example is also explained in 
Example \ref{bex}.
\end{example}

\medskip

 {\it Quantum integrable models} : certain Grothendieck rings can be identified with the ring of commuting operators of a quantum integrable model, via the transfer-matrix construction. Indeed, a quantum integrable model involves a space $W$ (the quantum space) with a family of commuting operators on this space (including quantum Hamiltonians). The transfer-matrix construction defines a ring morphism
$$T: K_0(\mathcal{C})\rightarrow \text{End}(W)$$
where $\mathcal{C}$ is the category of representation of an algebra of symmetry of the system. Hence, any relation in $K_0(\mathcal{C})$ (a "universal relation") gives relations between the eigenvalues of the operators on $V$, for any related quantum integrable model. 

\begin{example} The famous $XXZ$-model can be realized using the representation theory of the quantum affine algebra associated to $sl_2$ (see below).\end{example}

The aim of this lecture is to discuss certain symmetries of Grothendieck rings, and their consequences. The main examples of symmetries we present are Cluster symmetries and Weyl group symmetries.

\section{Cluster algebra symmetries}

We discuss the first symmetries of Grothendieck rings in this lecture, in the context of cluster algebras.

\subsection{A quick review on cluster algebras}

 The theory of cluster algebras was introduced by Fomin-Zelevinsky \cite{FZ}. A cluster algebra is a commutative algebra 
with a distinguished set of generators grouped into overlapping subsets (the clusters, obtained by an inductive process called mutation). Each element of the cluster algebra 
can be expressed as a rational fraction into the elements of a cluster. From a geometric point of view, one might 
think about a variety with various sets of local coordinates, and each element of the algebra of functions over 
this variety can be expressed locally in terms of the local coordinates. Then the mutations can be understood as a regular change of coordinates.

More precisely, the cluster algebra $\mathcal{A}_Q$ is a commutative algebra 
associated with a quiver $Q$ (without loops or 2-cycles). It has a distinguished set of generators (the cluster variables) defined by a combinatorial process.
Cluster algebras have many incarnations in various fields, and in particular in representation theory as we will discuss below.

The cluster algebra $\mathcal{A}_Q$ attached to the quiver
$Q$ (with set of vertices $Q_0$) is a subring of the field 
$$\mathcal{F} = \mathbb{Q}(X_i)_{i\in Q_0}.$$
with free variables $X_i$ which are called the initial cluster variables (together with the initial quiver $Q$, they form the initial cluster). The cluster algebra 
$\mathcal{A}_Q$ is defined as the subalgebra of $\mathcal{F}$ generated by the cluster variables, 
obtained inductively from the initial variables by an inductive process called mutations. 
For example, the first step mutated variables $X_i^*$ are defined by the formula : 
$$X_i X_i^* = \prod_{j\rightarrow i}X_j  + \prod_{j\leftarrow i}X_j,$$
were the arrows $\rightarrow$ and $\leftarrow$ are arrows in the initial quiver $Q$ (the quiver gets also mutated in the process). 
The number of cluster variables is not necessarily finite.
In addition, the cluster variables are grouped into overlapping subsets called clusters, which are all in bijection with $Q_0$. 
The cluster monomials are defined as the monomials into the cluster variables of the same cluster. In some situations, one 
may have additional non-mutable cluster variables : they are called frozen variables and they belong to all clusters.

One of the fundamental properties of cluster algebras is the Laurent phenomenon : any cluster variable can be expressed as a Laurent polynomial
in the cluster variables of a given seed. In the initial seed, this can be written as 
$$\mathcal{A}_Q\subset \mathbb{Z}[X_i^{\pm 1}]_{i\in Q_0}.$$

\begin{example} Consider a $Q$ quiver of type $A_2$ : 
$$\bullet \longrightarrow \bullet$$
We have the initial cluster variables $(X_1,X_2)$, and five cluster variables : 
$$X_1 , X_2, \frac{1 + X_2}{X_1}, \frac{1 + X_1}{X_2}  
, \frac{1 + X_1 + X_2}{X_1X_2}.$$
There are also five clusters : 
$$(X_1,X_2), \left(\frac{1 + X_2}{X_1}, X_2\right), 
\left(\frac{1 + X_2}{X_1}, \frac{1 + X_1 + X_2}{X_1X_2}\right),$$
$$\left(\frac{1 + X_1}{X_2},\frac{1 + X_1 + X_2}{X_1X_2}\right),  \left(\frac{1 + X_1}{X_2}, X_1\right).$$
\end{example}

\begin{rem} This example should not be misleading : in general, the number of clusters is not equal to the number of cluster variables.
\end{rem}

\subsection{Monoidal categorification of cluster algebras}

In this lecture we focus on the relation between cluster algebras and monoidal categories 
(there are also important relations with additive categories, see \cite{K} and references therein).

\begin{defi}\cite{HL1} A monoidal category $\mathcal{M}$ is a {\it monoidal categorification} of a cluster algebra $\mathcal{A}$ is there exists a ring isomorphism 
$$\phi : \mathcal{A}\rightarrow K(\mathcal{M})$$
so that the cluster monomials are sent to simple classes.
\end{defi}

In \cite{HL1}, a stronger version is also discussed, where it is required that $\phi$ induces a bijection between the cluster monomials and 
the real simple classes (whose tensor square is also simple). In the definition above, the simple modules corresponding to 
cluster monomials are called reachable modules. They have the property to be real, as the square of a cluster monomial is a cluster monomial.

A monoidal categorification can be seen as a cluster symmetry of the Grothendieck ring, or as a categorical realization of the cluster algebra.
It gives useful informations on the monoidal category : for example, it points out remarkable simple representations (corresponding to cluster variables) and it gives a factorization of simple representations corresponding to cluster monomials (in the same way as cluster monomials factorize as a product of cluster variables).

The original examples \cite{HL1} were obtained from the finite-dimensional representations of quantum affine algebras.
Many developments about monomial categorifications in this context have followed, one can cite for instance 
\cite{bc, KKKO, Nak:cluster, FQin}, and different contexts \cite{cw} for instance. 

In all known examples, the cluster symmetry is related to the defect of symmetry of the category which produces relations in the Grothendieck ring : we have objects $V$, $W$ so that 
$$V\otimes W\text{ not isomorphic to }W\otimes V$$
but 
$$[V\otimes W] = [W\otimes V]$$ 
in the Grothendieck ring (examples will be given below).

\subsection{Quantum affine algebras}

Consider $\mathfrak{g}$ a complex simple Lie algebra (for example, $\mathfrak{g} = \mathfrak{sl}_2$) of rank $n$ 
(that is $n = 1$ for $\mathfrak{g} = \mathfrak{sl}_2$). 
Then we can construct the affine Kac-Moody algebra $\hat{\mathfrak{g}}$, which is a 
central extension of the loop algebra\footnote{We do not consider the derivation element in this lecture.} 
$$\mathfrak{g}\otimes \mathbb{C}[t,t^{-1}].$$

For $q\in\mathbb{C}^*$ a quantum parameter (which is assumed to not root of unity), we have the corresponding quantum affine algebra : $\mathcal{U}_q(\hat{\mathfrak{g}})$. It is a Hopf algebra and a $q$-deformation of the universal enveloping algebra $\mathcal{U}(\hat{\mathfrak{g}})$.

Let $\mathcal{C}$ be the category of finite-dimensional representations of $\mathcal{U}_q(\hat{\mathfrak{g}})$ : 
it has a very rich structure. This monoidal category if not semi-simple and not braided.
The simple objects of $\mathcal{C}$ have been parametrized in terms of Drinfeld polynomials by Chari-Pressley, 
that is of $n$-tuples $(P_i(z))_{1\leq i\leq n}$ of monic polynomials. 
In particular, for $1\leq i\leq n$ and $a\in\mathbb{C}^*$, we have the fundamental representation $V_i(a)$ which corresponds to 
the $n$-tuple $(1,\cdots, 1, 1 - za, 1, \cdots, 1)$ with a degree $1$ polynomial in position $i$.

\begin{theorem}\cite{Fre} The Grothendieck ring 
$K(\mathcal{C})$ is commutative and polynomial 
$$K(\mathcal{C})\simeq \mathbb{Z}[X_{i,a}]_{1\leq i\leq n, a\in\mathbb{C}^*},$$
with $X_{i,a} = [V_{i,a}]$ class of a fundamental representation.
\end{theorem}

\subsection{Monoidal categorification - examples}\label{mce}

As a first example, consider the complex Lie algebra $\mathfrak{g} = \mathfrak{sl}_2$. We define $\mathcal{M}$ as the 
monoidal Serre subcategory of  the category $\mathcal{C}$ of finite-dimensional representations of $\mathcal{U}_q(\hat{\mathfrak{sl}}_2)$-representations generated by two fundamental representations
$$V_1(1), V_1(q^2).$$
Then we obtain a monoidal categorification
$$K(\mathcal{M})\simeq \mathcal{A}_Q$$
with $Q$ of type $A_2$ (as above) with one frozen variable. The total number of cluster variables is $3$. 
Indeed, there is a simple representation $W$ in $\mathcal{M}$ of dimension $3$ so that we have exact sequences 
$$0\rightarrow {\bf 1}\rightarrow V_1(q^2)\otimes V_1(1)\rightarrow W\rightarrow 0,$$
$$0\rightarrow W\rightarrow V_1(1)\otimes V_1(q^2)\rightarrow {\bf 1}\rightarrow 0,$$
which come from degenerated braidings (which are not isomorphisms as the two tensor products are not isomorphic) : 
$$V_1(q^2)\otimes V_1(1)\rightarrow  V_1(1)\otimes V_1(q^2)\text{ and }V_1(1)\otimes V_1(q^2)\rightarrow V_1(q^2)\otimes V_1(1).$$
These exact sequences correspond to the unique mutation relation in the cluster algebra $\mathcal{A}_Q$ : 
$$[V_1(q^2)][V_1(1)] = 1 + [W].$$
The two clusters are $(V_1(q^2), W)$ and $(V_1(1), W)$. As an application, we obtain that every simple representations in $\mathcal{M}$ can be factorized into these $3$ representations. Note that there are infinitely many simple classes in this category $\mathcal{M}$, which have the following form with $a,b\geq 0$ : 
$$(V_1(1))^{\otimes a}\otimes W^b\text{ and }(V_1(q^2))^{\otimes a}\otimes W^b.$$

As another example, consider the complex Lie algebra $\mathfrak{g} = \mathfrak{sl}_3$. We define $\mathcal{M}$ as the 
monoidal Serre subcategory of  the category $\mathcal{C}$ of finite-dimensional $\mathcal{U}_q(\hat{\mathfrak{sl}}_3)$-representations generated by the four fundamental representations
$$V_1(1), V_1(q^2), V_2(q), V_2(q^3).$$
Then we obtain a monoidal categorification
$$K(\mathcal{M})\simeq \mathcal{A}_Q$$
with $Q$ of type $A_2$ (as above) with two additional frozen variables. The total number of cluster variables is $7$.

As an application, we obtain that every simple representation in $\mathcal{M}$ can be factorized into these $7$ representations. 
Note that there are infinitely many simple classes in this category $\mathcal{M}$.

\subsection{Global structures}

For a general simple Lie algebra $\mathfrak{g}$, a cluster algebra structure has been obtained on $K(\mathcal{C}^-)$ for $\mathcal{C}^-$ a large subcategory\footnote{This large subcategory is sufficient to describe all simple modules by using spectral parameter shifts.} of $\mathcal{C}$ \cite{HL2} . Then Kashiwara-Kim-Oh-Park established the following conjectured in \cite{HL2}.

\begin{theorem}\cite{KKOP} The category $\mathcal{C}^-$ of finite-dimensional representations of a quantum affine algebra $\mathcal{U}_q(\hat{\mathfrak{g}})$ is a monoidal categorification of a cluster algebra.
\end{theorem}

As an application, we could prove for simple representations corresponding to cluster monomials, a Conjecture that I formulated in 2004 \cite{H1} for non simply-laced types. Using cluster algebra methods, the result extends results obtained by Nakajima \cite{Nak:quiver} in simply laced-types.

\begin{theorem}\cite{FHOO} The dimension (and character) of simple representations corresponding to cluster monomials can be obtained from an algorithm (\`a la Kazhdan-Lusztig).
\end{theorem}

\begin{rem} The algorithm is partly based on the existence of a quantization of the Grothendieck rings, introduced in \cite{Nak:quiver, VV:qGro} for simply-laced types using quiver varieties (and generalized in \cite{H1} using a different method). The deformation is then identified with a natural quantization of the cluster algebra (the quantum cluster algebra). Recently, with R. Fujita \cite{FuH}, we propose a construction of monoidal Jantzen filtrations for monoidal categories, with certain degenerated braidings, in order to give a purely categorical construction of such quantum Grothendieck rings. This is based on distinguished filtrations by submodules :
$$M = M_0\supset M_1\supset \cdots \supset M_N = \{0\}$$
which lead to deformed classes : 
$$[M]_t = \sum_{r\geq 0} t^r [M_r/M_{r-1}],$$
and to a deformation $*$ of the product of the Grothendieck ring, by considering filtrations of tensor products.

As an illustration, in the first example of Section \ref{mce}, we have the filtrations
$$V_1(q^2)\otimes V_1(1)\supset {\bf 1}\supset \{0\}\text{ and }V_1(1)\otimes V_1(q^2)\supset W\supset \{0\}$$
which lead to the deformed products
$$[V_1(q^2)]*[V_1(1)] = t [W] + 1\text{ and }[V_1(1)]*[V_1(q^2)] = t^{-1} [W] + 1.$$
\end{rem}

\section{Shifted quantum affine algebras}

The next step is to handle shifted quantum groups. The shifted quantum affine algebras form a new class of quantum groups closely related to (quantized $K$-theoretical) Coulomb branches(in the sense of Braverman-Finkelberg-Nakajima \cite{BFN}). They were introduced by Finkelberg-Tsymbaliuk \cite{FT} in the study of these $K$-theoretical Coulomb branches.

Consider $\mathfrak{g}$ and $q$ as above. The algebra $\mathcal{U}_q^\mu(\hat{\mathfrak{g}})$ can be seen as a variation of the quantum affine algebra 
$\mathcal{U}_q(\hat{\mathfrak{g}})$ depending on a shift parameter : a coweight $\mu$ of the Lie algebra $\mathfrak{g}$ 
(for example in the $sl_2$-case, we can see $\mu$ as an integer).
For $\mu = 0$, $\mathcal{U}_q^0(\hat{\mathfrak{g}})$ is essentially $\mathcal{U}_q(\hat{\mathfrak{g}})$.
These algebras have a very interesting representation theory.

\begin{theorem}\cite{H} $\mathcal{U}_q^\mu(\hat{\mathfrak{g}})$ has a non-zero finite-dimensional representation if and only if $\mu$ is codominant.
\end{theorem}

In general, $\mathcal{U}_q^\mu(\hat{\mathfrak{g}})$ has an abelian category $\mathcal{O}^\mu$ of representations (which might be 
finite-dimensional or infinite-dimensional). Recall $n$ is the rank of $\mathfrak{g}$.

\begin{theorem}\cite{H} The simple objects in $\mathcal{O}^\mu$ are parametrized by $n$-tuples of rational fractions $\Psib = (\psi_i(z))_{1\leq i\leq n}$ regular at $0$, with the degree condition : 
$$\text{deg}(\psi_i(z)) = \alpha_i(\mu)$$
where $\alpha_i$ is the simple root of $\mathfrak{g}$ attached to the vertex $i$ of the Dynkin diagram.
\end{theorem}

\begin{example}  For $1\leq i\leq n$, $a\in\mathbb{C}^*$, consider $\Psib$ defined by $\psi_j(z) = 1$ if $j\neq i$ and
$$\psi_i(z) = q^{r_i}\frac{1 - zaq^{-r_i}}{1 - zaq^{r_i}},$$
with $r_i$ the length of the simple root $\alpha_i$. Then we recover the finite-dimensional fundamental representation $V_i(a)$ 
of the ordinary quantum affine algebra (the rational fractions $\psi_j(z)$ should not be confused with the Drinfeld polynomials). 
\end{example}

\begin{example}
 For $1\leq i\leq n$, $a\in\mathbb{C}^*$, consider 
$$\Psib = ((1 - za)^{\delta_{i,j}})_{1\leq j\leq n}.$$ 
The associated simple representation of $\mathcal{U}_q^{\omega_i^\vee}(\hat{\mathfrak{g}})$ is called the positive prefundamental representation $L_{i,a}^+$ (here $\omega_i^\vee$ is a fundamental coweight). It is of dimension $1$ !  
\end{example}

\begin{example}\label{exbom} For $1\leq i\leq n$, $a\in\mathbb{C}^*$, consider 
$$\Psib = ((1 - za)^{-\delta_{i,j}})_{1\leq j\leq n}.$$ 
The associated simple representation of $\mathcal{U}_q^{-\omega_i^\vee}(\hat{\mathfrak{g}})$ is called the negative prefundamental representation $L_{i,a}^-$. It is infinite dimensional. In addition, $\mathcal{U}_q^{-\omega_i^\vee}(\hat{\mathfrak{g}})$ contains a copy of the ordinary quantum affine Borel algebra $\mathcal{U}_q(\hat{\mathfrak{b}})$. Then, restricted to this subalgebra, $L_{i,a}^-$ is the simple prefundamental representation constructed in \cite{HJ} as a limit of simple finite-dimensional representations. 
\end{example}

Without loss of generality, we may assume that all roots and poles of the rational fractions are integral powers of $q$.

One obtains an abelian category
$$\mathcal{O}^{sh} = \bigoplus_{\mu\text{ coweight}} \mathcal{O}^\mu$$
The (completed) sum of Grothendieck groups
$$K_0(\mathcal{O}^{sh}) = \bigoplus_{\mu} K_0(\mathcal{O}^\mu)$$
has a (topological) ring structure\footnote{It will be proved in \cite{HZ} that the subcategory of finite length representations 
is stable by fusion product, and so it gives rise to a (non-topological) ring structure on $K_0(\mathcal{O}^{sh})$.} induced by a the fusion product construction 
(which defines a representation $V* W$ in $\mathcal{O}^{\mu + \lambda}$ from 
simple representations $V$, $W$ in $\mathcal{O}^\mu$ and $\mathcal{O}^\lambda$ respectively).
This ring can be considered as an analogue of a Grothendieck ring. In addition, it is commutative.

Consider the subcategory $\mathcal{C}^{sh}\subset \mathcal{O}^{sh}$ of finite-dimensional representations of shifted quantum affine algebras. 
Its Grothendieck ring $K(\mathcal{C}^{sh})$ contains as a subring the Grothendieck ring $K(\mathcal{C})$ of finite-dimensional representations of the ordinary quantum affine algebra.

\begin{theorem}\cite{HL2, KKOP, H} $K_0(\mathcal{C}^{sh})$ is isomorphic to a cluster algebra $\mathcal{A}_{\Gamma_\infty}$, with an explicit quiver $\Gamma_\infty$. The initial cluster variables are classes of positive
prefundamental representations. The cluster monomials correspond to certain classes of simple modules.
\end{theorem}

\begin{example} For $\mathfrak{g} = \mathfrak{sl_2}$, the quiver $\Gamma_\infty$ is the infinite linear quiver :
$$\cdots\longrightarrow \bullet \longrightarrow \bullet \longrightarrow \cdots$$
The initial seed is formed of classes of $1$-dimensional positive prefundamental representations
$$\cdots\longrightarrow  L_{1,q^{-2}}^+ \longrightarrow L_{1,1}^+ \longrightarrow L_{1,q^{2}}^+ \longrightarrow \cdots$$
(up to factors by invertible representations, that we omit in these notes for clarity).
The first step mutations are given by certain remarkable relations called Baxter $TQ$-relations
$$[L_{1,a}^+][V_1(a)] = [L_{1,aq^2}^+] + [L_{1,aq^{-2}}^+]$$
where $V_1(a)$ : $2$-dimensional fundamental representation. In fact, this relation is also relevant from the point of view of quantum integrable models.
\end{example}

\begin{theorem}\cite{GHL} The Grothendieck ring $K_0(\mathcal{O}^{sh})$ is isomorphic to (a completion of) a cluster algebra $\mathcal{A}_{\Gamma_\infty'}$, with an explicit quiver $\Gamma_\infty'$.
\end{theorem}

 The proof is partly based on additional symmetries of the Grothendieck ring that will be discussed in the next section.

\begin{example} For $\mathfrak{g} = \mathfrak{sl}_2$, the new quiver $\Gamma_\infty'$ is obtained from the older quiver $\Gamma_\infty$ by just inverting the direction of one arrow : 
$$\cdots\longrightarrow \bullet \longrightarrow \bullet  {\color{red}\longleftarrow} \bullet \longrightarrow \bullet \longrightarrow \cdots$$
Then the initial seed is made of positive and negative prefundamental representations : 
$$\cdots\longrightarrow  L_{1,q^2}^+ \longrightarrow L_{1,1}^+ \longleftarrow L_{1,q^{-2}}^- \longrightarrow L_{1,q^{-4}}^-\longrightarrow \cdots$$
The mutation at the vertex corresponding to $L_{1,1}^+ $ 
$$\cdots\longrightarrow  L_{1,q^2}^+ \longleftarrow L_{1,1}^- \longrightarrow L_{1,q^{-2}}^- \longrightarrow L_{1,q^{-4}}^-\longrightarrow \cdots$$
as we have the remarkable $QQ$-relation 
\begin{equation}\label{remqq}[L_{1,1}^+][L_{1,1}^-] = 1 + [L_{1,q^2}^+][L_{1,q^{-2}}^+].\end{equation}
Similarly, the mutation at the vertex corresponding to $L_{1,q^{-2}}^-$ translated the seed in the other direction : 
$$\cdots\longrightarrow  L_{1,q^2}^+ \longrightarrow L_{1,1}^+ \longrightarrow L_{1,q^{-2}}^+ \longleftarrow L_{1,q^{-4}}^-\longrightarrow \cdots$$
\end{example}

In general, we conjecture the following. 

\begin{conj}\cite{GHL} All cluster monomials in $\mathcal{A}_{\Gamma_\infty'}$ correspond to classes of simple objects in $\mathcal{O}$ through our isomorphism.
\end{conj}

This general Conjecture is, for moment, only known in the following case.

\begin{theorem}\cite{GHL} The conjecture is true for $\mathfrak{g} = sl_2$.
\end{theorem}

\section{Symmetry of $q$-characters}

We now discuss certain group symmetries of Grothendieck rings.

\subsection{Classical Theory}

Let us briefly the well-known Weyl group symmetry of classical characters. Let $\mathfrak{g}$ be a complex finite-dimensional simple Lie algebra of rank $n$, and $K_0(\mathfrak{g})$ be the Grothendieck ring of its finite-dimensional representations. 
The character morphism, which encodes the dimension of weight spaces of finite-dimensional representations, defines an injective ring morphism 
$$\chi : K_0(\mathfrak{g})\rightarrow \mathbb{Z}[y_i^{\pm 1}]_{1\leq i\leq n}.$$
The image of $\chi$ can be characterized as an invariant subring
$$\text{Im}(\chi) = (\mathbb{Z}[y_i^{\pm 1}]_{1\leq i\leq n})^W$$
where the Weyl group $W$ is generated by the simple reflexions $s_i$, $1\leq i\leq n$ which act by the explicit formula
$$s_i(y_j) = y_j a_i^{-\delta_{ij}}\text{ where }
a_i = \prod_{1\leq k\leq n}y_k^{C_{ji}}.$$
Here $a_i$ corresponds to a simple root, and the $C_{ji}$ are coefficients of the Cartan matrix of $\mathfrak{g}$.

\begin{example} For $\mathfrak{g} = sl_2$, the Cartan matrix $C = (2)$ and we have
$$a_1 = y_1^2$$
$$s_1(y_1) = y_1 a_1^{-1} = y_1^{-1},$$
$$s_1^2(y_1) = y_1$$
$$s_1(y_1 + y_1^{-1}) = y_1 + y_1^{-1}.$$
$$\text{Im}(\chi) = (\mathbb{Z}[y_1^{\pm 1}])^W = \mathbb{Z}[ y_1 + y_1^{-1}].$$
\end{example}

\subsection{$q$-characters}

The analogue of characters for finite-dimensional representations of quantum affine algebras $\mathcal{U}_q(\hat{\mathfrak{g}})$ are $q$-characters (as defined by Frenkel-Reshetikhin \cite{Fre}). It defines an injective ring morphism on the Grothendieck ring $K(\mathcal{C})$ of the category $\mathcal{C}$ 
of finite-dimensional representations (here we assume that the quantization parameter $q\in\mathbb{C}^*$ is not a root 
of unity) : 
$$\chi_q : K(\mathcal{C}) \rightarrow \mathcal{Y} = \mathbb{Z}[Y_{i,a}^{\pm 1}]_{1\leq i\leq n, a\in\mathbb{C}^*}.$$
It is a refinement of characters for the action of the underlying finite-type quantum group $\mathcal{U}_q(\mathfrak{g})$ (as 
we can recover the ordinary by removing the spectral parameters $a$).

\begin{example} For $a\in\mathbb{C}^*$, the fundamental representation $V_1(a)$ of $\mathcal{U}_q(\hat{\mathfrak{sl_2}})$ satisfies : 
$$\chi_q(V_1(a)) = Y_{1,a} + Y_{1,aq^2}^{-1}.$$
If we forget the spectral characters, we recover the ordinary character $Y_1 + Y_1^{-1}$.
\end{example}

There are monomial analogs of simple roots (we give here the formula for simply-laced types) : 
$$A_{i,a} = Y_{i,aq^{-1}}Y_{i,aq}\prod_{j|C_{j,i} = -1} Y_{j,a}^{-1}.$$
For example for $\mathfrak{g} = \mathfrak{sl}_2$, we have $A_{1,a} = Y_{1,aq^{-1}}Y_{1,aq}$.

We may wonder what is the analogue of the Weyl group symmetry in this context.

\subsection{Symmetry of $q$-characters}

We introduce in \cite{FH2} the following operators : 
$$\Theta_i(Y_{j,a}) = Y_{j,a}A_{i,aq^{-1}}^{-\delta_{i,j}}\frac{\Sigma_{i,aq^{-3}}^{\delta_{i,j}}}{\Sigma_{i,aq^{-1}}^{\delta_{i,j}}},$$
for $1\leq i,j\leq n$ and $a\in\mathbb{C}^*$. Here $\Sigma_{i,a}$ is the solution of the $q$-difference equation
$$\Sigma_{i,a} = 1 + A_{i,a}^{-1}\Sigma_{i,aq^{-2}}$$
in a sum 
$$\Pi = \bigoplus_{w\in W}\tilde{\mathcal{Y}}^w$$ 
of completions $\tilde{\mathcal{Y}}^w$ of $\mathcal{Y}$, parametrized by the Weyl group.

\begin{example} For $\mathfrak{g} = sl_2$, we have
$$\Theta_1(Y_{1,a}) = Y_{1,aq^{-2}}^{-1}\frac{\Sigma_{1,aq^{-3}}}{\Sigma_{1,aq^{-1}}}.$$
Here $\Sigma_{1,a}$ is the couple
$$(1 + A_{1,a}^{-1}(1 + A_{1,aq^{-2}}^{-1}(1 + \cdots),$$ 
$$- A_{1,aq^2}(1 + A_{1,aq^4}(1 + \cdots))).$$
It belongs to 
$$\Pi = \mathcal{Y}^e \oplus \mathcal{Y}^{s_1}.$$
\end{example}

\begin{rem} One component of these operators is related to operators defined in \cite{In} on a different space.
\end{rem}

Now $\mathcal{Y}$ embeds in $\Pi$ diagonally, and we have the following analogue of the classical result.

\begin{theorem}\cite{FH2} The $\Theta_i$ are involutions and define a Weyl group action (with the simple reflexions acting as $\Theta_i$). Moreover, we have the ring of invariants
$$\mathcal{Y}^W = \text{Im}(\chi_q).$$
\end{theorem}

\subsection{Interpretation of the Weyl group action}

Let us explain an interpretation of the Weyl group action above, which connected to the cluster algebra symmetry.

\begin{theorem}\cite{FH, H} The replacement of the variables :  
$$Y_{i,a} = \frac{[L_{i,aq^{-1}}^+]}{[L_{i,aq}^+]}$$
in the $q$-character of finite-dimensional representation 
gives a well-defined relation in (the fraction field) of  $K(\mathcal{O}^{sh})$.
\end{theorem}

Note that originally, this Theorem was established in \cite{FH} in the context of representations of Borel subalgebra 
of quantum affine algebras (as in Example \ref{exbom}, in order to prove Frenkel-Reshetikhin conjecture on spectra of quantum integrable models), 
but it can also be stated for representations of shifted quantum affine algebras.

Motivated by quantum integrable models and Bethe Ansatz equations, we expect that the class $[L_{i,a}^+]$ can by replaced in the 
Theorem by the class of other simple representations, that we denote by $\tilde{L}_{i,a}^+$ (see \cite{FH31}). This is exactly how one can write the defining formula of the operator $\Theta_i$, which is now interpreted as a substitution of simple classes in the Grothendieck ring.

We can now relate this picture to cluster algebra symmetry. Indeed, starting from $Q_a^{\omega_i^\vee} = [L_{i,a}^+]$, and using the Weyl group action, one can introduce for a general $w\in W$, $i\in I$ :
$$Q_a^{w(\omega_i^\vee)} \in K(\mathcal{O}^{sh}).$$

\begin{theorem}\cite{FH3}[Generalized $QQ$-systems]
The series $Q_a^{w(\omega_i^\vee)}$ satisfy the $QQ$-system 
$$Q_{aq}^{(ws_i)(\omega_i^\vee)}Q_{aq^{-1}}^{w(\omega_i^\vee)} - Q_{aq^{-1}}^{(ws_i)(\omega_i^\vee)}Q_{aq}^{w(\omega_i^\vee)} 
 = \prod_{j|C_{i,j} = - 1} Q_a^{w(\omega_j^\vee)}.$$
\end{theorem}

Hence the same relation is satisfied by elements in an orbit under the Weyl group action.

\begin{example} For $\mathfrak{g} = sl_2$, this $QQ$-system is the quantum Wronskian relation : 
\begin{equation}\label{remqqs}\tilde{Q}_{aq}Q_{aq^{-1}} - \tilde{Q}_{aq^{-1}} Q_{aq} 
 = 1,\end{equation}
where we have denoted $Q_a = Q_a^{\omega_1^\vee}$ and $\tilde{Q}_a = Q_a^{-\omega_1^\vee}$. 
\end{example}

These $QQ$-system relations are crucial in the construction of the cluster algebra structure on $K(\mathcal{O}^{sh})$ in \cite{GHL}, 
as they are appear in mutations of the cluster algebra (for example, compare Equation (\ref{remqq}) and Equation (\ref{remqqs})). 
Hence, cluster symmetries are related to group symmetry of Grothendieck rings in these examples, which one could find quite intriguing.

\end{document}